# A New Implicit-Explicit Local Method to Capture Stiff Behavior with COVID-19 Outbreak Application

Huseyin Tunc, Murat Sari[1]

Department of Mathematics, Faculty of Arts and Science, Yildiz Technical University, Istanbul 34220, Turkey.

**Abstract**

In this paper, a new implicit-explicit local method with an arbitrary order is produced for stiff initial value problems. Here, a general method for one-step time integrations has been created, considering a direction free approach for integrations leading to a numerical method with parameter-based stability preservation. Adaptive procedures depending on the problem types for the current method are explained with the help of local error estimates to minimize the computational cost. Priority error analysis of the current method is made, and order conditions are presented in terms of direction parameters. Stability analysis of the method is performed for both scalar equations and systems of differential equations. The currently produced parameter-based method has been proven to provide $A-$stability, for $0.5 \leq \theta < 1$, in various orders. The present method has been shown to be a very good option for addressing a wide range of initial value problems through numerical experiments. It can be seen as a significant contribution that the Susceptible-Exposed-Infected-Recovered equation system parameterized for the COVID-19 pandemic has been integrated with the present method and stability properties of the method have been tested on this stiff model and significant results are produced. Some challenging stiff behaviours represented by the nonlinear Duffing equation, Robertson chemical system, and van der Pol equation have also been integrated, and the results revealed that the current algorithm produces much more reliable results than numerical techniques in the literature.

---

[1]Corresponding author: Tel: +90 212 383 43 60, E-mail: sarim@yildiz.edu.tr (Murat Sari)

# 1. Introduction

Capturing numerical behaviour of differential equations encountered in various fields of science including physics, chemistry, engineering, biology etc. is vital. In general, exact solutions of the equations cannot be derived easily all the time, or analytical expressions are very complicated to observe behaviour of the physical system. Even for the large linear equation systems, analytical evaluations are not easy to implement, and symbolic calculations lead to enormous computational time. These drawbacks of analytical approaches can be handled by considering accurate and economic numerical methods. In numerical techniques, it is in general easy to convert differential equations to algebraic equations and solve explicit or implicit recursive relations. By appropriate selection of parameters used in a numerical method, the convergence of the method can generally be controlled. However stiff differential equations are not easy task to handle and to find out reliable numerical solutions in the entire domain. Stiff behaviours are modelled by not only ODEs [1-2] but also PDEs [3-4] and arise in various phenomena. Numerical techniques are also classified according to the stiff and non-stiff problems in the literature [2,5]. The methods constructed to solve stiff problems are event of special interest due to the instability of classical method. Unwanted oscillations, divergence or slow convergence are commonly faced in solving stiff problems under the consideration of some inappropriate numerical techniques. Additionally, stiff ODE systems do not only come out directly in applications but also are encountered in discretization of some challenging processes described by PDEs [6].

In recent years, various numerical techniques have been derived in the literature to deal with challenging stiff problems [7-15]. A quasi-consistent fixed-step size numerical integration was used as an implicit stiff IVP solver in literature [7]. The method was found to be stable and was seen to preserve accuracy but is not computationally optimal due to the use of fixed step sizes. Another numerical study on solving stiff IVPs was presented as a one-leg explicit-implicit method [8], which is second order A-stable. Although this method offers an important numerical approach for stiff IVPs, the low degree of the method and assuming fixed step sizes lead to a computationally ineffective algorithm. For solving stiff differential-algebraic equations, an adaptive implicit Euler scheme was proposed in literature [9]. Even if the method produces acceptable results for stiff differential-algebraic equations, the order of the method is one, and the adaptive procedure requires a lot of time steps to achieve acceptable results over long intervals of time. Abdi et al. proposed arbitrary sequential barycentric rational finite difference method regarding the backward difference formulations (BDFs) [10]. This method class was shown to be adaptable to stiff IVPs, but the A-stable cases are obtained for at most second order formulations as it is also the case in the BDFs. A multistep method based on the BDFs has been derived for linear stiff problems in the literature [11]. The proposed method class offers a different perspective on the construction of constrained optimization in terms of the coefficients of BDFs rather than using predetermined coefficients. The authors tested their method on the linear equation, and the method seems to have a lack of computational efficiency for nonlinear stiff problems due to the extra cost of the optimization process. It is known from the literature that the discontinuous Galerkin methods (DGM) class is very effective for the differential equations representing shock behaviours. In the study of Fortin and Yakoubi [12], advantages of the DGM have been discussed for stiff IVPs and an adaptive formulation has been presented. A considerable disadvantage of the DGM is the coupling the degree of freedom at common interior nodes and this makes the method a bit computationally ineffective. Also, a modified approach based on radial basis functions that improve the quality of known numerical solvers was offered by Gu and Jung in [13] for non-stiff problems. In addition, recently some techniques such as various versions of Runge-Kutta methods [14-18], backward differentiation methods [19-21] and collocation methods [22-24] have been presented to solve stiff problems. All the above methods used to solve stiff IVPs in the literature are artificial because they use the differential equation taken into account in the derived formulae and, unlike our method, cannot get all the information from the differential equation.

As described in the literature above, a large number of continuous and discrete methods have been analysed and successfully implemented to deal with the problems of interest. Despite all these underlined advantages, those methods still have suffered from certain drawbacks such as having specific formulae as in finite difference-based algorithms, having enormous degrees of freedom, as in finite element-based techniques, having poor stability properties, as in differential transform-based methods. The present method overcomes the corresponding disadvantages with its arbitrarily high order and minimized degrees of freedom. As evidenced in the following sections, the current approach leads to optimization of the known information in each iteration and allows easy computation of adaptive meshes. This method, which we call the implicit-explicit local differential transformation method (IELDTM), is a multi-derivative time integration. In the proposed method, which is a stability-preserved numerical solver that accepts the idea of differential transformation as a starting point, all information about the numerical process is taken directly from the corresponding differential equations. An efficient method has thus here been produced for establishing parameter-dependent and adaptive implicit-explicit LDTM, which is appropriate for both stiff and non-stiff differential equations with optimum cost. It has been proved than the present numerical method offers much more better convergence properties than existing the semi-analytic DTM [25] and the explicit local DTM [26-28]. Priori error analysis of the present method has been done thoroughly, and order conditions are also produced. Stability analysis of the current method has been done for single equations and system of equations. Stability results are discussed by determining $A$−stable and $L$−stable states depending on direction parameters and it has been proved that the IELDTM contains up to fourth order $A$−stable schemes. Stability regions are demonstrated for special directions, i.e. for forward, central and backward cases. Depending on problem types, we have produced adaptive procedures for the current method to minimize computational effort. In this approach, we have used the known local differential transform values of the functions that come out in local error estimates to find a reliable adaptive procedure. It has been proved that the current numerical method has ability to accurately solve stiff IVPs without any considerable limitation. In numerical experiments, we have focused on various types of challenging stiff differential equations and illustrated the efficiency of the method. Since the stiff ODE solvers of MATLAB, *ode15s, ode23s, ode113 and ode45*, are widely and effectively used in the related problems, the effectiveness of the present method has been compared with those solving challenging stiff problems represented by the SEIR epidemiological system, the Van der Pol equation, the Robertson chemical system and the cubic nonlinear Duffing equation. In finding the results, the IELDTM has been observed to be far more cost-effective and accurate than the MATLAB stiff solvers, *ode15s, ode23s, ode113 and ode45*.

## 2. Implicit-explicit Local Differential Transform Method

The background of the local differential transformation can be seen in literature [26-28]. For simplicity, we omit the differential transform properties here. In this section, the IELDTM is introduced for a first order IVP with priori error analysis and stability analysis. Following the analysis of the currently produced method, adaptive procedures are presented for various problem types.

Consider the following IVP,

$$\boldsymbol{x}'(t) = G(\boldsymbol{x}(t), t), t > 0, \quad \boldsymbol{x}(0) = \boldsymbol{C}, \tag{1}$$

where $\boldsymbol{x}(t) \in R^m, \boldsymbol{C} \in R^m, G: R^m \times R \to R^m$ and the function $G(\boldsymbol{x}(t), t)$ satisfies the Lipschitz condition with constant $L$. We assume that the IVP (1) is well-posed and the exact solution $\boldsymbol{x}(t)$ is analytic on the considered domain.

Let us divide the interval $[0, T]$ into at most $N$ time elements with $\Delta t_i = t_{i+1} - t_i$ and the partition of the interval as $\omega = 0 = t_0 < t_1 < \cdots < N^* = T$ where $N^* \leq N$. Let us consider the convergent Taylor series representation of the function $x(t)$ about $t_i$ as

$$x_i(t) = \sum_{k=0}^{K} X_i(k)(t - t_i)^k + O((t - t_i)^{K+1}), \quad t_i - \rho^i \leq t \leq t_i + \rho^i \tag{2}$$

where $i = 0, 1, \ldots, N^*$ and $\rho^i$ is the radius of convergence of the representation. From our assumption, the function $x_i(t)$ is analytic and has the radius of convergence satisfying $\rho^i > \Delta t_i$. The function $x(t)$ has also a convergent Taylor series representation about $t_{i+1}$ with the radius of convergence at least $\rho^{i+1} = \rho^i - \Delta t_i$ and on the interval of convergence $(t_{i+1} - \rho^{i+1}, t_{i+1} + \rho^{i+1})$. Then, the function $x(t)$ can be written as

$$x_{i+1}(t) = \sum_{k=0}^{K} X_{i+1}(k)(t - t_{i+1})^k + O((t - t_{i+1})^{K+1}), \quad t_{i+1} - \rho^{i+1} \leq t \leq t_{i+1} + \rho^{i+1}. \tag{3}$$

Assuming $\rho^i > \Delta t_i$ and $\rho^{i+1} > \Delta t_i$, we can conclude that

- because of the convergence assumptions, two convergent representations (2) and (3) need to give the same numerical result at any point of the interval $[t_i, t_{i+1}]$.
- any point in the interval $[t_i, t_{i+1}]$ can be written as $t^* = t_i + (1 - \theta)\Delta t_i$, where $0 \leq \theta \leq 1$.
- the explicit-implicit method can be produced with the $C^0$-continuity of the solutions at such interior points.

Convergent solutions $x_i(t)$ and $x_{i+1}(t)$ at the interior points $t^* = t_i + (1 - \theta)\Delta t_i$ need to satisfy the following continuity relation

$$x_{i+1}(t_i + (1 - \theta)\Delta t_i) = x_i(t_i + (1 - \theta)\Delta t_i). \tag{4}$$

With the use of the representations, the following equality needs to be hold

$$\sum_{k=0}^{K} X_{i+1}(k)(-\theta \Delta t_i)^k = \sum_{k=0}^{K} X_i(k)((1 - \theta)\Delta t_i)^k + O((\Delta t_i)^{K+1}, \theta) \tag{5}$$

where $O((\Delta t_i)^{K+1}, \theta)$ represents the dependency of the local truncation error to time increment, transformation order and direction parameter. It is time to determine the algebraic relations between $X_i(k)$ and $X_{i+1}(k)$ for all $k = 0, 1, \ldots, K$.

Taking differential transform of (1), the following relation is obtained

$$X_i(k + 1) = \frac{1}{k+1} F(X_i(k), t_i) \tag{6}$$

where $i = 0, 1, \ldots, N^*$, $k = 0, 1, \ldots, K - 1$ and $F$ is the transformed form of the function $G(x, t)$ at the local point $t = t_i$. It is obvious that $X_i(k)$ can be written in terms of $X_i(0)$ for each $i$ and $k$. Thus, using the algebraic relations between $X_i(k)$ and $X_{i+1}(k)$ and putting into (5) the explicit-implicit equation

$$g(X_{i+1}(0), \theta, \Delta t_i) = h(X_i(0), \theta, \Delta t_i) \tag{7}$$

is found. Here the functions $g$ and $h$ are obtained from the right and left hand sides of Eq. (5). In the right side of (7), all parameters $X_i(0), \theta$ and $\Delta t_i$ are known because of the previous step. The left side of (7) depends explicit/implicitly on the unknown $X_{i+1}(0)$. Thus, (7) is either linear or nonlinear equation(s) of $X_{i+1}(0)$ depending on the function $G(x, t)$.

Then setting the initial condition $x_0(0) = X_0(0) = c$ and predetermining adaptive parameters $\theta$ and $\Delta t_i$, (7) can be solved either directly or numerically for each step. Note that we assume the following global solutions

$$x(t_i) = X_i(0)$$

for each $i = 0, 1, \ldots, N^*$. Therefore, the following remarks should be underlined because of its importance for the rest of the study:

- $\theta = 0$ with $N^* = 1$ leads to the classical semi-analytic DTM [25].
- $\theta = 0$ with arbitrary number of $N^*$ leads to the explicit forward scheme and known as local or multi-step differential transform method [26-29].
- The rest of the selections of the parameter $\theta$, i.e. $\theta \neq 0$, gives rise to the currently derived implicit schemes. To the best of the authors' knowledge, this derivation has not been studied in the DTM.
- $\theta = 0.5$ yields the produced implicit central scheme with arbitrary order here, which gives us a stable and order preserved method for both stiff and non-stiff cases.
- $\theta = 1$ leads to a derivation of the implicit backward scheme which is both stable and order preserved method for both stiff and non-stiff cases.
- Stability preserving schemes can be obtained by considering $\theta \geq 0.5$.
- The classical $\theta$−method including Crank-Nicolson method is a special case of the IELDTM; $\theta = 0.5$ and $K = 1$.

In the following subsection, the priori error analysis of the derived numerical technique is done by illustrating the order conditions.

## 2.1 Error Analysis

Consider the following initial value problem,

$$x'(t) = G(x(t), t), t > 0, \quad x(0) = C, \tag{8}$$

where $C \in R^m$ and $G: R^m \times R \to R^m$. Considering the procedure explained in the previous section and assuming the fixed time increment $\Delta t_i = \Delta t$, the present scheme leads to,

$$X_n(k + 1) = \frac{1}{k+1}[F(X_n(k), t_n)], k = 0, 1, 2, \ldots, K - 1, n = 0, 1, 2, \ldots, N, \tag{9}$$

$$x(t_{n+1}) \cong x_{n+1} = X_{n+1}(0) = \sum_{k=0}^{K} X_n(k)\big((1-\theta)\Delta t\big)^k - \sum_{k=1}^{K} X_{n+1}(k)(-\theta\Delta t)^k \qquad (10)$$

where $X_n$ and $F$ are the transformed forms of the functions $x_n(t)$ and $G(x_n(t), t)$. The exact solution of (8) at point $t = t_{n+1}$ can be expressed in the Taylor expansion form as

$$x(t_{n+1}) = X_{n+1}(0) = \sum_{k=0}^{K} X_n(k)\big((1-\theta)\Delta t\big)^k - \sum_{k=1}^{K} X_{n+1}(k)(-\theta\Delta t)^k + \Delta t \rho_n \qquad (11)$$

where $n = 0, 1, 2, \ldots, N-1$. Local truncation error can be obtained from the residuals of the Taylor expansions as follows

$$\rho_n = [(1-\theta)^{K+1} X_n(K+1) - (-\theta)^{K+1} X_{n+1}(K+1)] \Delta t^K + O(\Delta t^{K+1}). \qquad (12)$$

Now using the expansion $X_{n+1}(K+1) = X_n(K+1) + \Delta t(K+2) X_n(K+2) + O(\Delta t^2)$ leads to

$$\rho_n = [(1-\theta)^{K+1} X_n(K+1) - (-\theta)^{K+1} X_n(K+1)] \Delta t^k + O(\Delta t^{K+1}). \qquad (13)$$

Thus, the present scheme is of order $K+1$ if $\theta = 1/2$ and $K$ is odd and is of order $K$ for other selections of parameter $\theta$ and $K$. To have convergent numerical scheme, we also need to have bounded global error of the scheme. For further analysis about global discretization error, assume that the problem is linear, $G(x(t), t) = Ax(t) + B(t)$ and let $\varepsilon_n = x(t_n) - x_n$ for $n = 0, 1, \ldots, N$. The present method for this equation yields the following discretization,

$$\sum_{k=1}^{K} X_{n+1}(k)(-\theta \Delta t)^k = \sum_{k=0}^{K} X_n(k)\big((1-\theta)\Delta t\big)^k. \qquad (14)$$

With the use of recursive relation (5) with $G(x(t), t) = Ax(t) + B(t)$, the general term $X_n(k)$ can be stated as

$$X_n(k) = \frac{1}{k!} A^k X_n(0) + \sum_{p=0}^{k-1} \frac{p!}{k!} F(t_n, p) \qquad (15)$$

where $F$ is transformed form of the function $B(t)$. Then, substituting (15) into (14) gives the following recursive form

$$\sum_{k=0}^{K} \frac{1}{k!} A^k (-\theta\Delta t)^k X_{n+1}(0) = \sum_{k=0}^{K} \frac{1}{k!} A^k \big((1-\theta)\Delta t\big)^k X_n(0) + \sum_{k=0}^{K} \sum_{p=0}^{k-1} \frac{p!}{k!} \Big[ F(p, t_n)\big((1-\theta)\Delta t\big)^k - F(p, t_{n+1})(-\theta\Delta t)^k \Big]. \qquad (16)$$

Let us define the following stability functions

$$R_1(\Delta t A, \theta) = I + (1-\theta)\Delta t A + \big((1-\theta)\Delta t\big)^2 A^2 + \cdots + \big((1-\theta)\Delta t\big)^K A^K \qquad (17)$$

$$R_2(\Delta t A, \theta) = I - \theta \Delta t A + (\theta \Delta t)^2 A^2 + \cdots + (-\theta \Delta t)^K A^K. \qquad (18)$$

With the use of (15) into also exact expansion (11) and subtraction of (16) from (11) gives the following relation

$$\varepsilon_{n+1} = R(\Delta t A, \theta) \varepsilon_n + \delta_n. \qquad (19)$$

where $\delta_n = (R_2)^{-1}\Delta t \rho_n$ and $R(\Delta t A, \theta) = (R_2)^{-1} R_1$. Here the important fact is that the scheme is of order $p$ if $||\delta_n|| = O(\Delta t^{p+1})$. To bound global discretization error in terms of initial error $\varepsilon_0$ and local discretization error, recursive relation (19) becomes,

$$\varepsilon_n = R^n \varepsilon_0 + \sum_{p=0}^{n-1} R^{n-p-1} \delta_p. \tag{20}$$

Stability of the scheme is also depending on the bound of the norm estimate $||R^n||$ for all $n\Delta t \leq T$, where $T$ is the final time. So, assuming the following stability criteria,

$$||R^n|| \leq S, \text{ for all } n \geq 0 \text{ and } n\Delta t \leq T, \tag{21}$$

$$||(R_2)^{-1}|| \leq P, \tag{22}$$

the error norm bound becomes

$$||\varepsilon_n|| \leq S||\varepsilon_0|| + S \sum_{p=0}^{n-1} ||\delta_p||. \tag{23}$$

Then finally, using the definition of local discretization error $||\delta_p|| \leq P\omega \Delta t^{K+1}$ for all $p$, $\omega$ is defined as

$$\omega = \begin{cases} \left| \frac{1}{(K+2)!} \max_{t \in [0,T]} \left( \left| \frac{dx^{K+2}(x(t),t)}{dt^{K+2}} \right| \right) \left( \frac{1}{2} \right)^{K+1} (K+1) \right|, \text{if } \theta = 1/2 \text{ and } K \text{ is odd} \\ \left| \frac{1}{(K+1)!} \max_{t \in [0,T]} \left( \left| \frac{d^{K+1}x}{dt^{K+1}} \right| \right) [(1-\theta)^K - (-\theta)^K ] \right|, \text{if } \theta \neq \frac{1}{2} \text{ and } K \text{ is arbitrary.} \end{cases}$$

Assumption of $x(0) = x_0$ leads to the following error norm estimates

$$||x(t_n) - x_n|| \leq \omega^* \Delta t^{K+1}, \text{ if } \theta = 1/2 \text{ and } K \text{ is odd} \tag{24}$$

$$||x(t_n) - x_n|| \leq \omega^* \Delta t^K, \text{ otherwise} \tag{25}$$

where $\omega^* = SP\omega t_n$ for all $n = 1, 2 ..., N$. Thus, whenever the exact solution is smooth and stability criteria (21)-(22) are satisfied, then the present IELDTM converges to the exact solution with order $K + 1$ or $K$ depending on the selection of the parameters $\theta$ and $K$.

## 2.2 Stability Analysis

In order to visualize the behaviour of the currently presented IELDTM, we have analysed the stability properties over the scalar test equation and the system of linear equations.

### 2.2.1 Scalar Test Equation

At first, let us consider the following scalar-complex test equation

$$x'(t) = \lambda x(t), x(0) = x_0 \tag{26}$$

with $\lambda \in C$. Application of the IELDTM, with order $K$, to this test problem gives

$$x_{n+1} = R(\Delta t \lambda) x_n,$$

$$R(z) = \frac{1+(1-\theta)z+\frac{((1-\theta)z)^2}{2!}+\cdots+\frac{((1-\theta)z)^K}{K!}}{1-\theta z+\frac{(\theta z)^2}{2!}+\cdots+\frac{(-\theta z)^K}{K!}} \cong e^z, \ z \to 0 \tag{27}$$

If we perturb the initial value $x_0$ to $\tilde{x}_0$, we get the recursion $\tilde{x}_{n+1} = R(\Delta t \lambda)\tilde{x}_n$ then the difference $\tilde{x}_n - x_n$ leads to the following stability requirement

$$\tilde{x}_n - x_n = R(\Delta t \lambda)^n (\tilde{x}_0 - x_0) \tag{28}$$

and therefore, stability region of the present method is the set

$$S = \{z \in C : |R(z)| \le 1 \}$$

where $S \subset C$. To analyse the current approach, we need to remind the maximum modulus theorem [1]. In the light of this theorem, we have concluded that:

- The central schemes, $\theta = 1/2$ of order $K = 1, 2, 3, 4$, are $A$-stable.
- The backward schemes, $\theta = 1$ of order $K = 1, 2$ are $A$-stable as well as $L$-stable.
- The schemes, with $\theta < 1/2$ of any order including LDTM ($\theta = 0$), have conditional stability.
- For higher order implicit schemes, $\theta \ge 1/2$, we have almost $A$-stability that means just the stability condition fails in a very little region in $\mathbb{C}^-$.

The stability regions of the currently produced method, the IELDTM, have been illustrated with $\theta = 0, \theta = 1/2$ and $\theta = 1$ for various orders, respectively, as seen in Figs. 1-3. As realized in the Figures, stability properties of the central and backward cases are much more stronger than the forward cases.

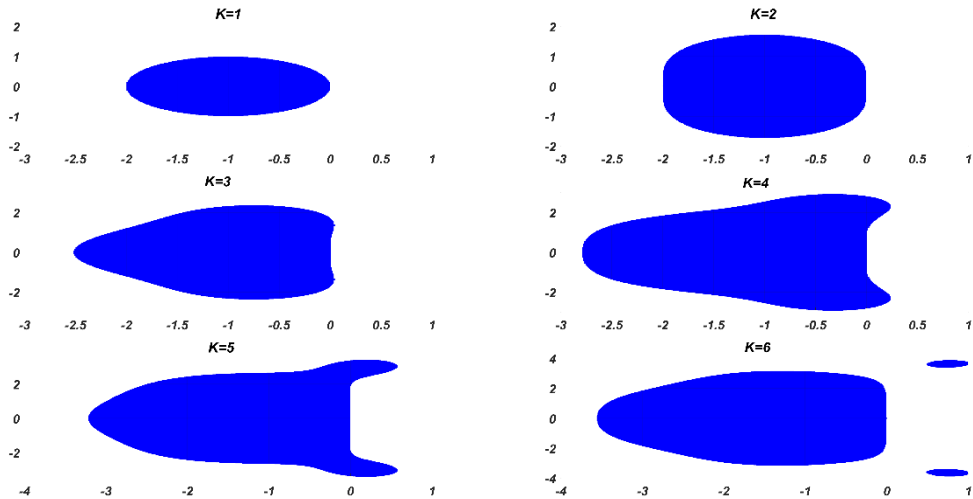

**Fig. 1**. Stability region (blue) of the forward schemes, $\theta = 0$, for various orders on the complex plane.

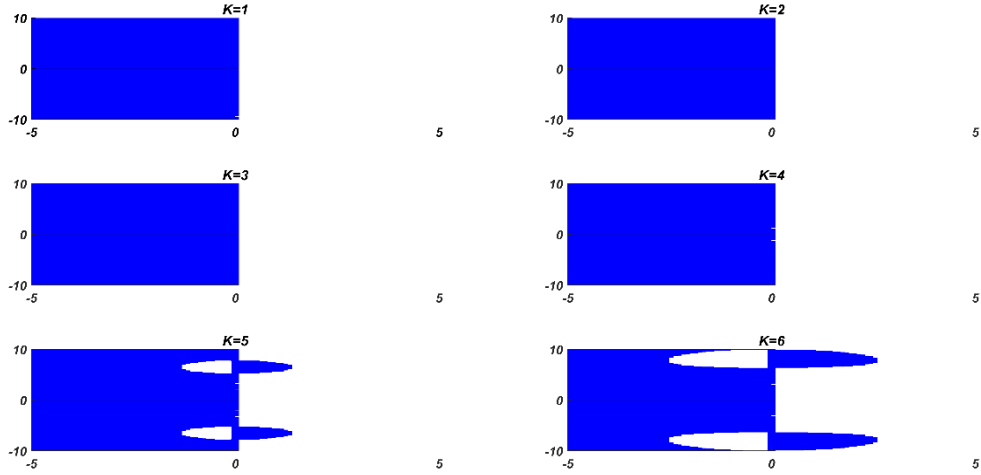

**Fig. 2.** Stability region (blue) of the central schemes, $\theta = 0.5$, for various orders on the complex plane.

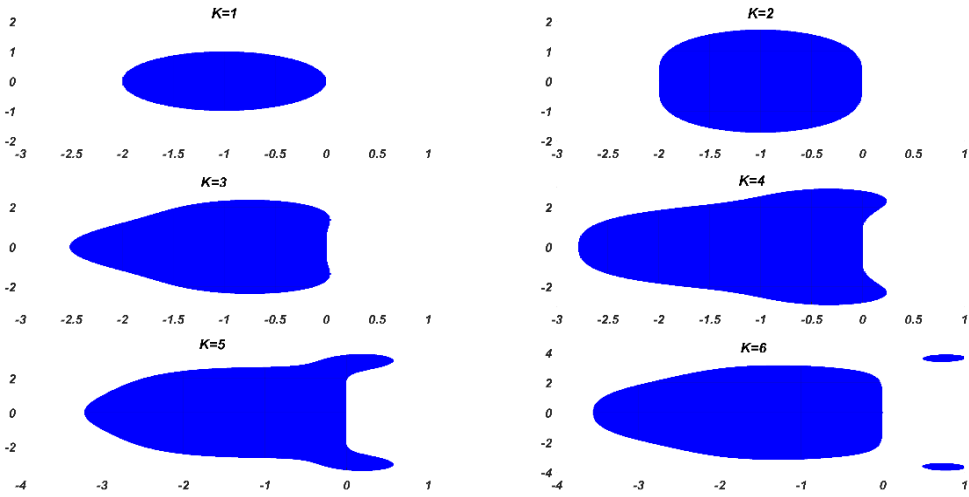

**Fig. 3.** Stability region (blue) of the backward schemes, $\theta = 1$, for various orders on the complex plane.

### 2.2.2 Stability for Linear Systems

Consider the following linear homogeneous system

$$\boldsymbol{x}'(t) = A\boldsymbol{x}(t), \ \boldsymbol{x}(0) = \boldsymbol{x}_0 \tag{29}$$

with $A \in R^{m \times m}$ and $\boldsymbol{x}_0 \in R^m$. Application of the IELDTM yields the following localized differential transformation coefficients

$$X_n(k) = \frac{1}{k!} A^k X_n(0). \tag{30}$$

Writing these coefficients into (5) leads to

$$\sum_{k=0}^{K}\left[\frac{1}{k!}A^{k}X_{n+1}(0)\right](-\theta\Delta t_{i})^{k} = \sum_{k=0}^{K}\left[\frac{1}{k!}A^{k}X_{n}(0)\right]\left((1-\theta)\Delta t_{i}\right)^{k} \tag{31}$$

where $X_n$ is the transformed form of $x(t)$. We define the stability functions as

$$R(\theta,\Delta tA) = \left[I - \theta\Delta tA + \frac{(\theta\Delta tA)^2}{2} + \cdots + \frac{(-\theta\Delta tA)^K}{K!}\right]^{-1}\left[I + (1-\theta)\Delta tA + \frac{((1-\theta)\Delta tA)^2}{2} + \cdots + \frac{((1-\theta)\Delta tA)^K}{K!}\right]. \tag{32}$$

Recalling $x(t_{n+1}) = X_{n+1}(0)$ and $x(t_n) = X_n(0)$, then (31) becomes

$$x(t_{n+1}) = R(\theta,\Delta tA)x(t_n). \tag{33}$$

For a perturbed system with initial condition $\bar{x}_0$, we get the same formula with (33). Subtraction of (33) from the corresponding perturbed system leads to

$$x(t_n) - \bar{x}(t_n) = R(\theta,\Delta tA)^n(x_0 - \bar{x}_0). \tag{34}$$

Thus, the stability function $R(\theta,\Delta tA)^n$ determines how much initial error will grow. Fortunately, we have a rich literature about stability conditions of like Eqn. (32) as all details analysed in [2]. First, we assume the following important theorem related to the bound of $||R||$ [2].

**Theorem 1.** Let the rational function $R(z)$ be bounded for $Re(z) \leq 0$ and assume that the matrix $A$ satisfies

$$Re(y, Ay) \leq 0 \text{ for all } y \in \mathbb{C}^n \tag{35}$$

where $(.,.)$ denotes the Euclidean scalar product. Then, in the matrix norm corresponding to the scalar product we have

$$||R(A)|| \leq \underbrace{\sup}_{Re(z) \leq 0} |R(z)|. \tag{36}$$

Note that for the homogeneous system (29), the following diminishing property leads to condition (35) as follows

$$\frac{d}{dt}||x||^2 = \frac{d}{dt}(x,x) = 2Re(x,x') = 2Re(x,Ax). \tag{37}$$

When the condition (35) is not satisfied, depending on the logarithmic norm of the matrix $A$, the following Corollary becomes useful for the bound of stability function $R$ [2].

**Corollary 1.** If a matrix occurred in the system (29) satisfies $Re(v,Av) \leq s||v||^2$ for all $v \in \mathbb{C}^m$ with $\mu = \mu(A) \leq s$ is the logarithmic norm of the matrix $A$, then

$$||R(A)|| \leq \underbrace{\sup}_{Re(z) \leq s} |R(z)|. \tag{38}$$

**Proof.** Assumption of $\bar{R} = R(z+s)$ and $\bar{A} = A - sI$ for Theorem 1 leads to the desired result.

In the light of Corollary 1, the present $A$-stable methods will be unconditionally stable for linear systems (29) whenever the logarithmic norm of the matrix A satisfies $\mu(A) \leq 0$ as will be stated in the next Corollary.

**Corollary 2**. If $\mu(A) \leq 0$, the $A$-stable cases of the IELDTM are unconditionally stable for homogeneous system (29).

**Proof.** Application of Corollary 1 to Eq. (34) gives

$$||x(t_n) - \bar{x}(t_n)|| \leq \left( \sup_{Re(z) \leq 0} |R(z)| \right)^n ||x_0 - \bar{x}_0||. \tag{39}$$

The $A$-stability assumption leads to the following estimate

$$\sup_{Re(z) \leq 0} |R(z)| \leq 1, \text{ for all } z \in \mathbb{C}. \tag{40}$$

Thus, the $A$-stable cases of the IELDTM are seen to be unconditionally stable. The table below provides information on some cases of $A$-stability of the IELDTM.

Note that for the implicit schemes, $\theta \geq 0.5$, with the rest of the orders have not been mentioned in Table 1, are almost $A$-stable as illustrated in the scalar case. Except that the central or backward cases of the IELDTM with $\theta \geq 0.5$ stated in Table 1, it is possible to find out various $A$-stable cases up to fourth-order approximation.

**Table 1.** Some special $A$-stable cases of the IELDTM

| $\theta$ | K |
|---|---|
| 0.5 | 1 |
| 0.5 | 2 |
| 0.5 | 3 |
| 0.5 | 4 |
| 1 | 1 |
| 1 | 2 |

## 2.3 Adaptive Searching

The IELDTM needs to selection of suitable parameters $\Delta t_i$ and $\theta$ to control both accuracy and stability. This situation also depends on the nature of the considered physical model. The established numerical approach is seen to be general implicit-explicit method depending on the parameter $\theta$. Since the stability properties of the implicit schemes and computational costs have been explained, an optimized technique needs to be obtained in terms of accuracy, stability and computational cost. To do this, the following adaptive techniques are produced depending on the nature of the problem.

**Case 1:**

In this case we assume that the problem is nonlinear and non-stiff. In non-stiff problems, the current forward methods generally produce convergent and computationally effective results. Thus, assuming explicit schemes are preferable in the present method by selecting the parameter as $\theta = 0$. Even if the model equations are nonlinear, algebraic equations can directly be solved without using any symbolic nonlinear solver. An adaptive procedure is here proposed for the time increment $\Delta t_i$ to reduce the computational cost. Considering the transformation of order $K$ with the selection parameter $\theta = 0$ leads to global error of order $K$ as we proved in error analysis. For a given tolerance the following criteria need to be satisfied for scalar equations as analysed in the literature [30],

$$\Delta t_i < \left(\frac{tol}{|X_i(K+1)|}\right)^{\frac{1}{K}}. \tag{41}$$

where $tol$ is predetermined tolerance and $X_i(K+1)$ can be obtained from the recursive relation. For a system of equations, the relation (41) can be modified to the related matrices and vectors by considering suitable norms. In this study, the norm $||.||_\infty$ is considered for the selection of adaptive $\Delta t_i$.

**Case 2:**

In this case, stiffness of the problem is assumed to be arbitrary. In such cases $\theta = 0.5$ and $\theta = 1$ will lead to accurate, stable and effective approaches. It is obvious from the error analysis that the error bounds of implicit schemes, $\theta \neq 0$, depend on both $X_i(K)$ and $X_{i+1}(K)$. The error bounds of central schemes have been analysed and thus concluded that taking odd values of the order $K$ is preferable. Then similar analysis for the central schemes, $\theta = 0.5$, with odd orders leads to the following bounds

$$\Delta t_i < \left(\frac{tol}{\left|\left(\frac{1}{2}\right)^{K+1}(K+1)\right)X_i(K+2)\right|}\right)^{\frac{1}{K+1}}. \tag{42}$$

where $tol$ is the predetermined tolerance and $X_i(K+2)$ is obtained from the recursion. For the backward adaptive cases, $\theta = 1$, the estimate (41) is also valid for the arbitrary selection of the transformation order $K$. The scalar adaptive procedure given in equation (42) can be enlarged to vectors and matrices with the use of suitable norms.

## 3. Numerical Experiments

Considering the quantitative and qualitative results of various test problems, the numerical representation and efficiency of the derived method are presented here. To measure the effectiveness of the current method in terms of accuracy and stability, challenging physical processes represented by stiff differential equations such as the stiff SEIR equation system, cubic nonlinear Duffing equation, Robertson equation system, and van der Pol equation have been dealt with. The computed results are compared with the results produced with FDM [10], *ode23s*, *ode15s*, *ode45* and *ode113* [31] as well as the exact solutions. To evaluate the error norms of the current results, absolute pointwise errors $E_i$ and maximum error norms $\|E\|_\infty$ have been preferred. The trust-region dogleg method [32] based *fsolve* built-in function of the MATLAB program have been applied to solve the resulting nonlinear algebraic equation system.

**Problem 1 [33]** Consider the following Susceptible-Exposed-Infected-Recovered/Removed (SEIR) system of equations

$$\frac{dS}{dt} = -\beta \left(\frac{S(t)}{N}\right)[P(t) + D(t) + \mu A(t)]$$

$$\frac{dE}{dt} = \beta \left(\frac{S(t)}{N}\right)[P(t) + D(t) + \mu A(t)] - \frac{1}{d_1}E(t)$$

$$\frac{dP}{dt} = \alpha \frac{E(t)}{d_1} - \frac{P(t)}{d_2}$$

$$\frac{dA}{dt} = (1-\alpha)\frac{E(t)}{d_1} - \frac{A(t)}{d_3} \tag{43}$$

$$\frac{dD}{dt} = \frac{P(t)}{d_2} - \frac{D(t)}{p}$$

$$\frac{dR}{dt} = \frac{D(t)}{p} + \frac{A(t)}{d_3}$$

where $S(t) + E(t) + P(t) + A(t) + D(t) + R(t) = N$ is the total population, $S(t), E(t), P(t), A(t), D(t), R(t)$ denote susceptible, exposed, pre-symptomatic, asymptomatic, hospitalized and recovered/removed number of individuals in the population, respectively. The parameters $\beta$, $\mu$ and $\alpha$ are daily transmission rate, transmission reduction factor and pre-symptomatic ratio, respectively. $d_1$, $d_2$, $d_3$ and $p$ represent the mean latency period, mean pre-symptomatic infectiousness period, mean asymptomatic infectiousness period, and mean hospitalization period, respectively. The SEIR model (43) is one of the key compartment models describing the spread of any infectious disease, such as the COVID-19 pandemic [33], in a population. If the basic reproduction number $R_0 = \beta[\alpha(d_2 + p) + (1-\alpha)d_3]$ is greater than one, then the number of susceptible individuals will asymptotically decrease to a threshold value with a stiff or non-stiff dynamic depending on the transmission rate $\beta$.

Taking differential transform of equation system (43) write as follows:

$$X_i(k+1) = \frac{1}{k+1}\left(AX_i(k) + B_i(k)\right) \qquad (44)$$

where

$$X_i(k) = [S_i(k), E_i(k), P_i(k), A_i(k), D_i(k), R_i(k)]^T,$$

$$A = \begin{bmatrix} 0 & 0 & 0 & 0 & 0 & 0 \\ 0 & -\frac{1}{d_1} & 0 & 0 & 0 & 0 \\ 0 & \frac{\alpha}{d_1} & -\frac{1}{d_2} & 0 & 0 & 0 \\ 0 & \frac{1-\alpha}{d_1} & 0 & -\frac{1}{d_3} & 0 & 0 \\ 0 & 0 & \frac{1}{d_2} & 0 & -\frac{1}{p} & 0 \\ 0 & 0 & 0 & \frac{1}{d_3} & \frac{1}{p} & -\alpha \end{bmatrix}$$

and

$$B_i(k) = \begin{bmatrix} -\frac{\beta}{N}\sum_{j=0}^{k} S_i(k-j)[P_i(j) + D_i(j) + \mu A_i(j)] \\ \frac{\beta}{N}\sum_{j=0}^{k} S_i(k-j)[P_i(j) + D_i(j) + \mu A_i(j)] \\ 0 \\ 0 \\ 0 \\ 0 \end{bmatrix}$$

for all $k = 0,1,\ldots,K-1$ and $i = 0,1,\ldots,N-1$. The parameter values $d_1 = 3.69$, $d_2 = 3.47$, $d_3 = 3.47$, $p = 1.92$, $\alpha = 0.14$, $\beta = 1.12$ and $N = 3.10^6$, which are obtained by Li *et. al.* (2020) for the dynamics of COVID 19 pandemics, considered for all simulations. The SEIR model does not describe a stiff behavior for moderate values of the transmission rate $\beta$.

The advantages of the IELDTM (with $\theta = 0$) over the existing LDTM [26-27,34-35] can be seen in Figure 4 to solve the moderately-stiff SEIR model (43). The LDTM suffers from instability when rapid decreasing is observed in the susceptible population. As compared with the reference *ode15s* solution of the system, the behaviours of the state variables are accurately captured by the central IELDTM using the same parameter values. The IELDTM with $\theta \geq 0.5$ destroys the instability drawback of the LDTM by solving the problem implicitly.

To increase the stiffness of the SEIR model, the transmission rate is modified to have the following time-dependent form

$$\beta(t) = \begin{cases} 1.12, & t \leq t_c \\ 1.12\eta, & t > t_c \end{cases} \qquad (45)$$

where η ≥ 1 is the scaling factor of the transmission rate. If $t_c$ is selected in an interval in which rapid decreasing of the susceptible population occurs, then the scaling factor η increases the stiffness of the problem. The SEIR model is solved by the adaptive central IELDTM and various MATLAB *ode* solvers with $t_c = 66$ and η ∈ [1, 12] and the results are illustrated in Figure 5. The required number of time steps are compared with respect to the changing stiffness of the problem in Figure 5a. It has been observed that the sixth and eighth order central IELDTMs need much less time elements than the adaptive MATLAB *ode23s*, *ode15s*, *ode113*, and *ode45*. Even if the number of required time steps of the MATLAB solvers increases with increasing $\eta$ values, the IELDTM is not affected by increasing stiffness. The behavior of the stiff dynamics of the susceptible population is illustrated in Figure 5b. It is observed that the central IELDTM accurately captures the dynamics at optimum cost.

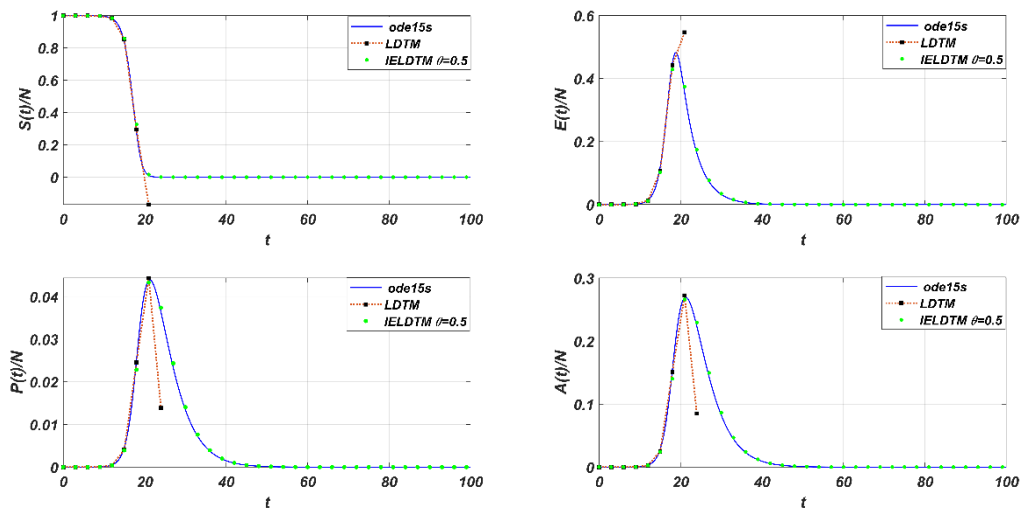

**Fig. 4.** Comparison of the central IELDTM and LDTM with the reference *ode15s* solutions of the SEIR model (43) for $K = 6$ and $\Delta t = 3$.

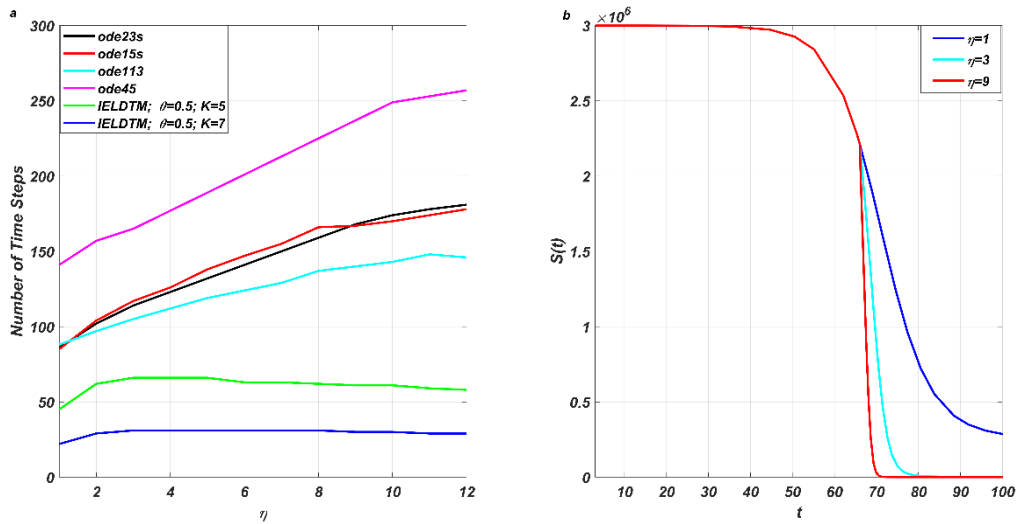

**Fig. 5.** a) Comparison of the required number of time steps of the central adaptive IELDTM, *ode23s*, *ode15s*, *ode113* and *ode45* for solving SEIR model (43) with the changing values of η and $tol = 10^{-5}$, b) The effect of the stiffness factor η on the dynamics of the susceptible population $S(t)$ is captured by the central adaptive IELDTM with an optimized degrees of freedom.

**Problem 2 [27]** Consider the following cubic nonlinear Duffing equation

$$x'' + \alpha x' + \beta x + \gamma x^3 = 0 \text{ for } 0 \leq t \leq t_f \tag{46}$$

with initial displacement and velocity

$$x(0) = x_0 \text{ and } x'(0) = x_0^*. \tag{47}$$

In the study of Tunc and Sari (2019b), it is proved that the logistic function $x(t) = \frac{1}{1+e^{-t}}$ is exact solution of Eq. (46) with the following parameters,

$$\alpha = -3, \ \beta = 2, \gamma = -2, \ x_0 = 0.5, \ x_0^* = 0.25 \tag{48}$$

With the use of $x_1(t) = x(t)$ and $x_2(t) = \dot{x}(t)$, the cubic nonlinear Duffing equation is transformed to the following nonlinear system

$$x_1'(t) = x_2(t) \tag{49}$$

$$x_2'(t) = -\alpha x_2(t) - \beta x_1(t) - \gamma x_1^3(t). \tag{50}$$

Taking differential transformation of (49)-(50) one can obtain

$$X_i(k+1) = \frac{1}{k+1}\big(AX_i(k) + B_i(k)\big) \tag{51}$$

Where $k = 0,1, \dots, K-1$, $X_i(k) = [(X_1)_i(k), (X_2)_i(k)]^T$ is the differential transformation of $[x_1(t), x_2(t)]^T$, $B_i(k) = \left[0, \gamma \sum_{l=0}^{k} \sum_{n=0}^{l} (X_2)_i(n)(X_2)_i(l-n)(X_2)_i(k-l)\right]^T$ and the constant matrix $A$ can be defined as follows

$$A = \begin{bmatrix} 0 & 1 \\ -\beta & -\alpha \end{bmatrix}. \tag{52}$$

The IELDTM yields the following equation

$$\sum_{k=0}^{K} X_{i+1}(k)(-\theta \Delta t_i)^k = \sum_{k=0}^{K} X_i(k)\big((1-\theta)\Delta t_i\big)^k \tag{53}$$

where $X_0(0) = [x_0, x_0^*]^T$, $i = 0,1, \dots, N-1$ and the coefficients can be evaluated from Eq. (51).

The performance of the current IELDTM is shown by comparing the maximum errors depending on the varying values of the direction parameter $\theta$ (see Figure 6). As theoretically expected, the maximum errors are reduced exponentially by increasing the order of the method. It can be observed from Figure 6 that choosing $\theta = 0.5$ with odd-numbered $K$ values leads to an $(K + 1) - th$ order numerical method. The present adaptive central scheme is used for various transformation orders from $K = 3$ to $K = 9$ with the estimate $tol = 10^{-15}$ and the numerical results are compared in terms of both accuracy and computational costs in Fig. 7. As seen in Fig. 7, increasing values of the parameter $K$ cause the computational cost to decrease as expected. The theoretical order expectations and experimental order averages of the current IELDTM are compared with various values of the transformation order $K$ and direction parameter $\theta$ as given in Table 2. There, the theoretical expectations are almost matching with the experiments. The current adaptive central scheme and the *ode15s-ode23s* solvers are compared with the maximum errors in Table 3 to measure the advantages of the IELDTM over the existing stiff *ode* solvers. In all comparisons, the central adaptive IELDTM appears to give better results with both high accuracy and using less time step than the literature [31].

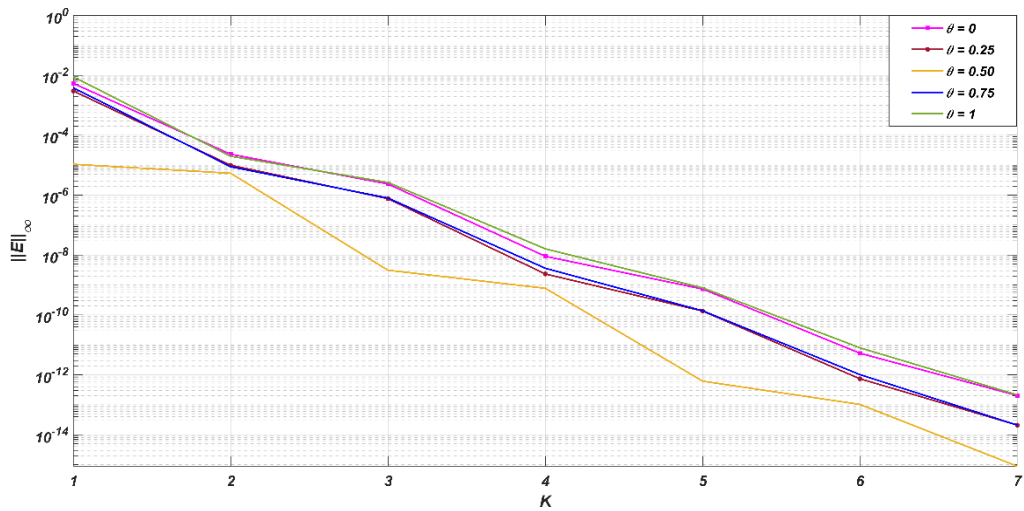

**Fig. 6.** Order refinement results of the IELDTM at $t_f = 1$ for various values of the direction parameter $\theta$ with the time increment $\Delta t = 0.05$.

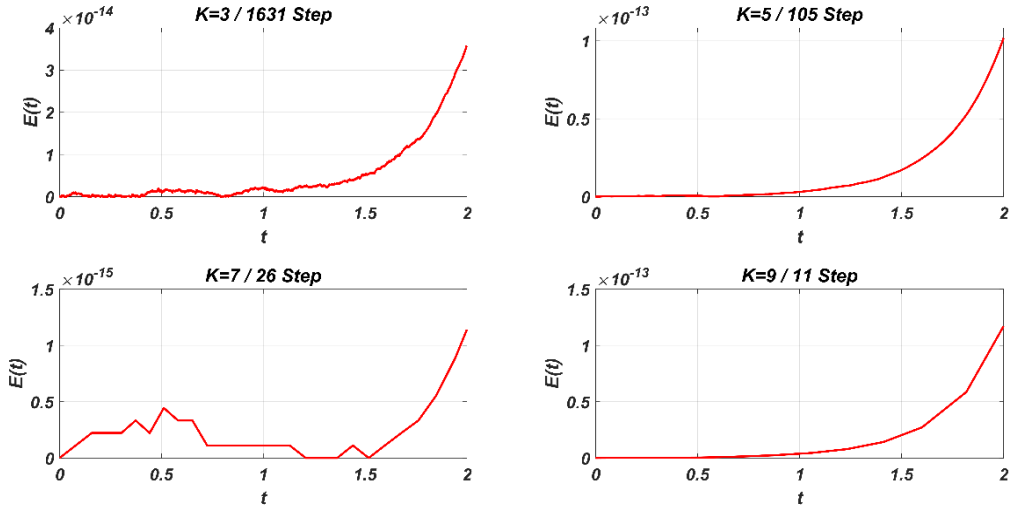

**Fig. 7.** Pointwise errors computed by the central adaptive IELDTM with various values of $K$ and $tol = 10^{-15}$

**Table 2** Effect of the approximation order of the IELDTM with various direction parameters in comparative way for Problem 2

| | $\theta = 0$ | | $\theta = 0.5$ | | $\theta = 1$ | |
|---|---|---|---|---|---|---|
| $K$ | Experimental Average Order | Theoretical Expectation | Experimental Average Order | Theoretical Expectation | Experimental Average Order | Theoretical Expectation |
| $K = 1$ | 0.8348 | 1 | 1.9804 | 2 | 1.2583 | 1 |
| $K = 2$ | 2.0603 | 2 | 1.9970 | 2 | 1.9319 | 2 |
| $K = 3$ | 2.9458 | 3 | 4.0142 | 4 | 3.0477 | 3 |
| $K = 4$ | 3.7102 | 4 | 4.0041 | 4 | 4.2008 | 4 |
| $K = 5$ | 4.9576 | 5 | 5.9965 | 6 | 5.0311 | 5 |
| $K = 6$ | 5.8559 | 6 | 6.0821 | 6 | 6.2854 | 6 |

**Table 3** Comparison of the present ACM-N and the *ode15s-ode23s* solvers with $tol = 10^{-10}$ for the adaptive time steps and maximum errors in Problem 2

| | IELDTM, $K = 3$ | | IELDTM, $K = 5$ | | *ode15s* [31] | | *ode23s* [31] | |
|---|---|---|---|---|---|---|---|---|
| $t_f$ | Step Number | Maximum Error | Step Number | Maximum Error | Step Number | Maximum Error | Step Number | Maximum Error |
| 1 | 51 | 7.93E-11 | 9 | 2.38E-10 | 64 | 1.37E-09 | 825 | 8.61E-08 |
| 2 | 92 | 4.62E-09 | 16 | 8.45E-09 | 99 | 3.08E-08 | 1487 | 2.99E-06 |
| 4 | 152 | 1.01E-05 | 27 | 1.49E-05 | 161 | 6.53E-05 | 2595 | 6.47E-03 |

**Problem 3 [10]** Consider the following modified Robertson chemical stiff system of nonlinear differential equations

$$x_1' = -0.04x_1 + 10^4 x_2 x_3 - 0.96 e^{-t},$$
$$x_2' = 0.04x_1 - 10^4 x_2 x_3 - 3.10^7 x_2^2 - 0.04 e^{-t}, \quad (54)$$
$$x_3' = 3.10^7 x_2^2 + e^{-t},$$

where $t \in [0, 4]$ and the initial conditions are

$$x_1(0) = 1, x_2(0) = x_3(0) = 0. \quad (55)$$

As stated in literature (Abdi et. al. 2019), the system (54) has the exact solution $x_1(t) = e^{-t}$, $x_2(t) = 0$ and $x_3(t) = 1 - e^{-t}$. Taking of differential transformation of (54) leads to the following relations,

$$(X_1)_i(k+1) = \frac{1}{k+1}\left[-0.04(X_1)_i + 10^4 \sum_{s=0}^{k}(X_2)_i(s)(X_3)_i(k-s) - 0.96\frac{e^{-t_i}}{k!}\right]$$

$$(X_2)_i(k+1) = \frac{1}{k+1}\Big[0.04(X_1)_i - 10^4 \sum_{s=0}^{k}(X_2)_i(s)(X_3)_i(k-s) -$$
$$3.10^7 \sum_{s=0}^{k}(X_2)_i(s)(X_2)_i(k-s) - 0.04\frac{e^{-t_i}}{k!}\Big] \quad (56)$$

$$(X_3)_i(k+1) = \frac{1}{k+1}\left[3.10^7 \sum_{s=0}^{k}(X_2)_i(s)(X_2)_i(k-s) + \frac{e^{-t_i}}{k!}\right]$$

where $(X_1)_0(0) = 1$, $(X_2)_0(0) = 0$, $(X_3)_0(0) = 0$, $k = 0, 1, \ldots, K-1$ and $i = 0, 1, \ldots, N-1$. Thus, the IELDTM yields the following equations

$$\sum_{k=0}^{K}(X_1)_{i+1}(k)(-\theta \Delta t_i)^k = \sum_{k=0}^{K}(X_1)_i(k)\big((1-\theta)\Delta t_i\big)^k$$
$$\sum_{k=0}^{K}(X_2)_{i+1}(k)(-\theta \Delta t_i)^k = \sum_{k=0}^{K}(X_2)_i(k)\big((1-\theta)\Delta t_i\big)^k \quad (57)$$
$$\sum_{k=0}^{K}(X_3)_{i+1}(k)(-\theta \Delta t_i)^k = \sum_{k=0}^{K}(X_3)_i(k)\big((1-\theta)\Delta t_i\big)^k$$

where $i = 0, 1, \ldots, N-1$ and the coefficients can be calculated from Eq. (56).

The present central adaptive schemes and the MATLAB solvers, *ode15s-ode23s* [31], have been compared to demonstrate the efficiency of the present approach through the pointwise errors (Figure 8). The current adaptive IELDTM schemes of order 4 and 6 have been shown to produce more accurate results, using far less number of time elements, than the commonly used stiff solvers, *ode15s-ode23s*. Maximum errors of the current ACM-N have been presented in comparison with the LBRFDM-BDF of the literature [10] (see Table 4). In the table, it is seen that the current schemes give more accurate results using the same number of time elements. It is also obvious from Table 4 that the IELDTM is far more accurate than various versions of the rational finite difference techniques presented for stiff behavior in the literature [10].

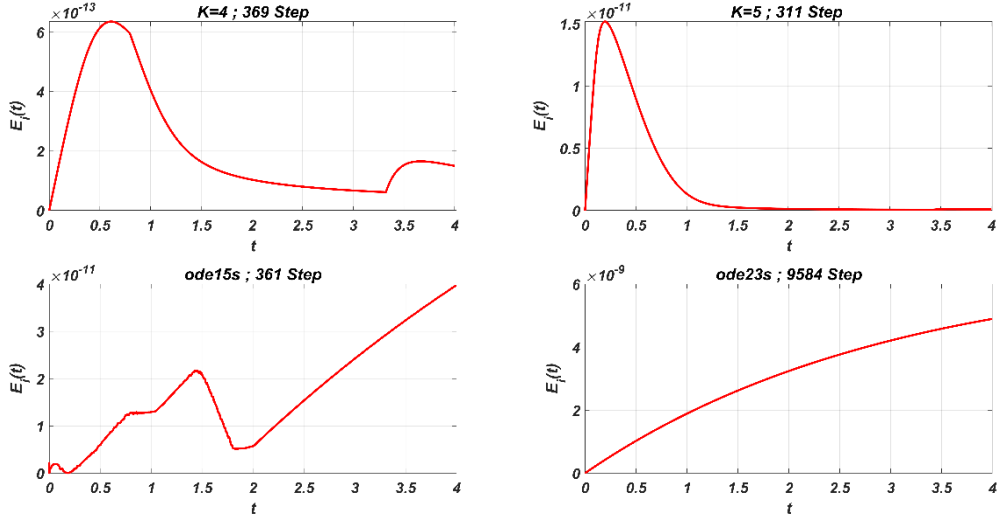

**Fig. 8.** Comparison of the present central adaptive IELDTM and the *ode15s-ode23s* through pointwise errors with $tol = 10^{-12}$ for Problem 3.

**Table 4** Comparison of the present central IELDTM with the literature LBRFDM-BDF [10] through maximum errors at $t = 4$ for Problem 3.

| $\Delta t$ | IELDTM $K = 3$ | IELDTM $K = 4$ | IELDTM $K = 5$ | LBRFDM (6,3) | LBRFDM (6,4) | 4-Step BDF |
|---|---|---|---|---|---|---|
| $2^{-5}$ | 2.69E-10 | 4.89E-11 | 3.89E-13 | 3.11E-07 | 1.55E-07 | 1.41E-07 |
| $2^{-6}$ | 4.97E-11 | 5.86E-12 | 3.79E-13 | 2.13E-08 | 1.07E-08 | 1.02E-08 |
| $2^{-7}$ | 4.97E-12 | 5.94E-13 | 1.33E-15 | 1.40E-09 | 6.98E-10 | 6.83E-10 |
| $2^{-8}$ | 4.76E-13 | 5.62E-14 | 8.88E-16 | 8.94E-11 | 4.47E-11 | 4.42E-07 |

**Problem 4.** Consider the Van der Pol equation in the first-order system form,

$$U' = V \tag{58}$$
$$V' = -U + \varepsilon(1 - U^2)V$$

where the initial conditions are taken to be as $U(0) = 2$ and $V(0) = 0$. Note that the system (58) corresponds to the second-order nonlinear Van der Pol equation with dependent variable $U$. The Van der Pol equation is a model of self-sustained oscillations of a triode electric circuit with the current voltage characteristics. Stiffness of the Van der Pol system (58) strongly depends on the parameter value $\varepsilon$. If the parameter $\varepsilon$ is not too large, then the equation is not stiff and can be solved with also non-stiff methods. However, in case of $\varepsilon \gg 1$, the system (58) is not easy to handle and requires more accurate numerical methods for its solution. Taking of differential transform of the system (58) yields

$$\boldsymbol{U}_i(k+1) = A\boldsymbol{U}_i(k) + \boldsymbol{F}(k) \tag{59}$$

where $\boldsymbol{U}_i(k) = [U_i(k), V_i(k)]^T$, $A = [0,1; -1, \varepsilon]$, $\boldsymbol{F}(k) = \left[0, -\varepsilon \sum_{l=0}^{k} \sum_{n=0}^{l} U_i(n) U_i(l-n) V_i(k-l) \right]^T$ and $k = 0, 1, \ldots, K-1$. The IELDTM leads to the following equations

$$\sum_{k=0}^{K} \boldsymbol{U}_{i+1}(k)(-\theta \Delta t_i)^k = \sum_{k=0}^{K} \boldsymbol{U}_i(k)\big((1-\theta)\Delta t_i\big)^k \tag{60}$$

where $i = 0, 1, \ldots, N-1$ and $\boldsymbol{U}_i(k)$ can be evaluated from the recursive relation (59).

non-stiff cases of the Van der Pol equation have been illustrated through the present central adaptive technique. The qualitative behavior is produced by considering the parameter values $K = 7$, $\theta = 0.5$, $t_f = 20$ and $tol = 10^{-10}$. As seen in Fig. 10, as the stiffness increases, the number of time steps required by the algorithm increases. The computed results of the stiff cases are presented as seen in Fig. 10 for the parameter values of $\varepsilon = 10$ and $\varepsilon = 100$, respectively, in Problem 4. The central adaptive technique is used for Fig. 10 with the parameter values $K = 5$, $\theta = 0.5$ and $tol = 10^{-10}$. As shown in the figures, with the current techniques, it is seen that the stiff behaviors are correct and successfully captured. The obtained results showed that the effectiveness of the present method IELDTM has been realized by comparing the MATLAB solvers, *ode15s-ode23s* [31] (see Table 5). In the table, it is understood that our method needs much less time step than the *ode* solvers and still appears to give more accurate results.

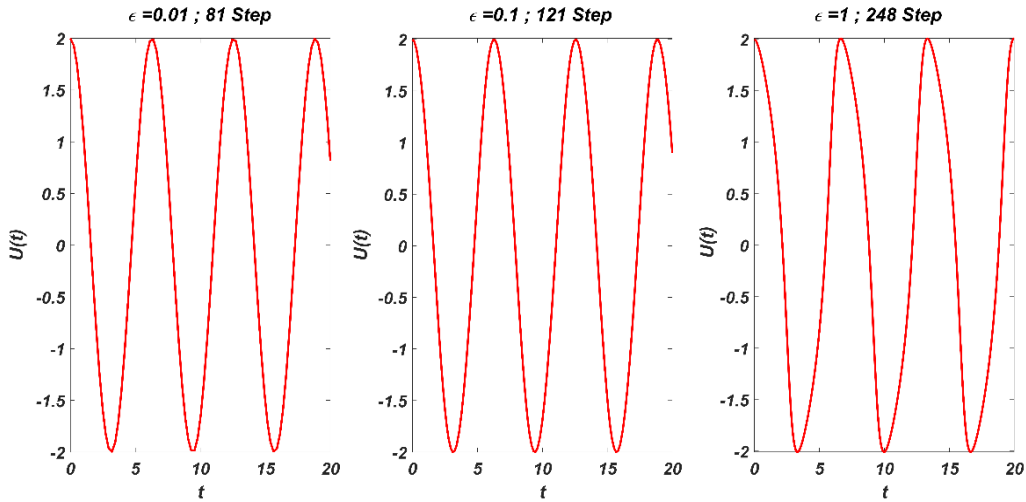

**Fig. 9.** Adaptive central solutions produced with $K = 7$, $\theta = 0.5$, $tol = 10^{-10}$ and various $\varepsilon$ values for Problem 4.

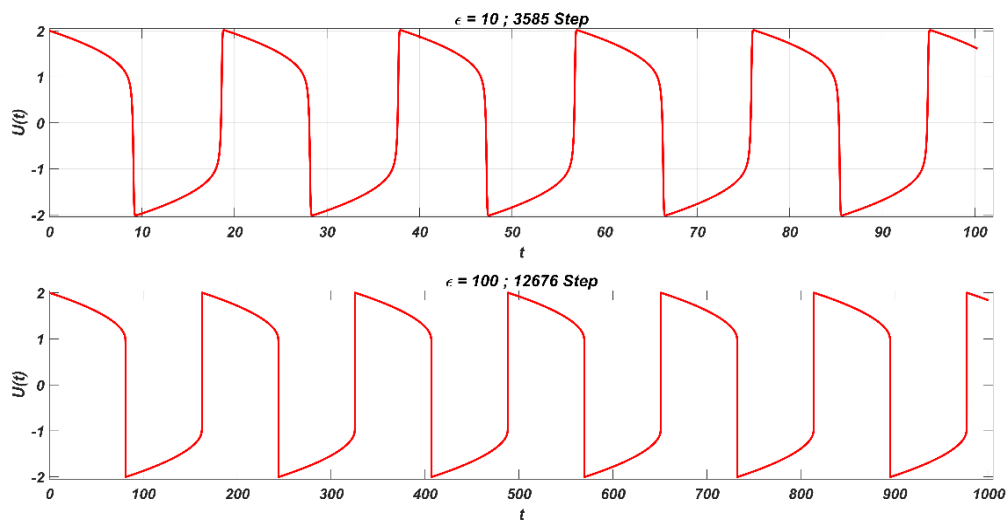

**Fig. 10.** Adaptive central solutions produced with $K = 5$, $\theta = 0.5$, $tol = 10^{-10}$, $\varepsilon = 10$ and $\varepsilon = 100$ for Problem 4.

**Table 5** Comparison of the present central adaptive IELDTM and the MATLAB solvers [31], *ode15s-ode23s*, with the adaptive time step numbers for $tol = 10^{-10}$ in Problem 4.

| $\varepsilon/T$ | IELDTM $K = 3$ | IELDTM $K = 5$ | IELDTM $K = 7$ | IELDTM $K = 9$ | *ode15s* | *ode23s* |
|---|---|---|---|---|---|---|
| 0.1/1 | 1788 | 254 | 95 | 53 | 700 | 10912 |
| 1/10 | 3827 | 520 | 193 | 108 | 1256 | 15717 |
| 10/100 | 45607 | 5339 | 1888 | 1068 | 9632 | 141653 |
| 100/1000 | 173179 | 19012 | 15282 | 10820 | 17523 | 299643 |

## 5. Conclusions

In this paper, a new one-step implicit-explicit local differential transform methods (IELDTM) have been developed with an arbitrary order for especially stiff initial value problems. The strong stability of the numerical method produced using the Taylor series has proven to be preserved. A priori error estimates of the currently derived approaches have been constructed and order conditions of the methods have been determined depending on direction parameters. Stability of the present methods has been discussed and they have thus been examined by considering $A-$ and $L-$stabilities. To reach optimal numerical methods, adaptive procedures have been produced for forward, central, and backward cases, respectively. The currently produced versatile methods have been seen to be effective for very challenging problems defined by stiff differential equations. To explain and analyse the challenging aspects of the problems, four striking stiff differential equations have been taken into account and thus qualitative and quantitative results have been exhibited. In a comparison way, efficiency of the present methods in terms of accuracy, stability, versatility and computational cost have been proved. In summary, the main contributions of the article are highlighted as follows:

- A numerically reliable, adaptive, high-order and stability preserved method has been derived,
- Error analysis and order conditions have been determined,
- Depending on the nature of the problem, various adaptive procedures have been proposed,
- $A$-stable and $L$-stable cases of the IELDTM have been identified, and the regions of stability have been found,
- Reliability of the adaptive procedures produced has been demonstrated and therefore embedding adaptivity features into the IELDTM by increasing the order of the method has been shown to minimize the computational cost.
- Versatility of the present algorithms in terms of challenging problem types has been discussed in detail.
- It has been observed that the SEIR equation system parameterized for the COVID-19 outbreak can be effectively integrated with the IELDTM and therefore the stiffness of the problem can be discussed in detail with this method

**Acknowledgments**

The first author would like to thank the Science Fellowships and Grant Programmes Department of TUBITAK (TUBITAK BIDEB) for their support to his academic research. The second author thanks Dr. Aniela Balacescu (Constantin Brâncuși University of Targu Jiu) for the great hospitality during the conference.

**References**

[1] W. Hundsdorfer, J. Verwer, Numerical solution of time-dependent advection-diffusion-reaction equations, Springer, 2003.

[2] E. Hairer, G. Wanner Solving ordinary differential equations 2- stiff and differential-algebraic equations, Springer, 1996.

[3] B. Kleefeld, J.M. Vaquero, SERK2v3: Solving mildly stiff nonlinear partial differential equation. J. Comput. Appl. Math. 299 (2016) pp.94-206.

[4 S. Li, Canonical Euler splitting method for nonlinear composite stiff evolution equations, Appl. Math. Comput. 289 (2016) pp.220-236.

[5] E. Hairer, S.P. Nørsett, G. Wanner, Solving ordinary differential equations 1 - nonstiff problems, Springer, 1996.

[6] M. Sari, H. Tunc, M. Seydaoglu, Higher order splitting approaches in analysis of the Burgers equation. Kuwait J. Sci. 46(1) (2019) pp.1-14.

[7] G.Y. Kulikov, R. Weiner, Doubly quasi-consistent fixed-stepsize numerical integration of stiff ordinary differential equations with implicit two-step peer methods, J. Comput. Appl. Math. 340 (2018) pp.256–275.

[8] Z Tan, C. Zhang, Implicit-explicit one-leg methods for nonlinear stiff neutral equations, Appl. Math. Comput. 335 (2018) pp.196–210.

[9] N. Leterrier, ARES: An efficient approach to adaptive time integration for stiff differential algebraic equations, Comput. Chem. Eng. 119 (2018) pp.46–54.

[10] A. Abdi, S.A. Hosseini, H. Podhaisky, Adaptive linear barycentric rational finite differences method for stiff ODEs, J. Comput. Appl. Math. 357 (2019) pp.204-214.

[11] B.V. Faleichik, Minimal residual multistep methods for large stiff non-autonomous linear problems, J. Comput. Appl. Math. 387 (2021) pp.112498.

[12] A. Fortin, D. Yakoubi, An adaptive discontinuous Galerkin method for very stiff systems of ordinary differential equations, J. Comput. Appl. Math. 358 (2019) pp.330–347.

[13] J. Gu, J.H. Jung Adaptive radial basis function methods for initial value problems, J. Sci. Comput. 82 (2020) pp.47.

[14] I. Higueras, T. Roldán, Strong stability preserving properties of composition Runge–Kutta schemes, J. Sci. Comput. 80 (2019) pp.784-807.

[15] L. Isherwood, Z.J. Grant, S. Gottlieb, Strong stability preserving integrating factor two-step Runge–Kutta methods, J. Sci. Comput. 81 (2019) pp.1446–1471.

[16] M. Narayanamurthi, Sandu, Efficient implementation of partitioned stiff exponential Runge-Kutta methods, Appl. Numer. Math. 152 (2020) pp.141-158.

[17] B.C. Vermeire, Paired explicit Runge-Kutta schemes for stiff systems of equations, J. Comput. Phys. 393 (2019) pp.465–483.

[18] C.A. Kennedya, M.H. Carpenter, Diagonally implicit Runge–Kutta methods for stiff ODEs, Appl. Numer. Math. 146 (2019) pp.221–244.

[19] Z.B. Ibrahim, K.I. Othman, M. Suleiman, Implicit r-point block backward differentiation formula for solving first-order stiff ODEs, Appl. Math. Comput. 186 (2007) pp.558–565.

[20] J.R. Cash, Modifed extended backward differentiation formulae for the numerical solution of stiff initial value problems in ODEs and DAEs, J. Comput. Appl. Math. 125 (2000) pp.117-130.

[21] O.A. Akinfenwa, S.N. Jator, N.M. Yao, Continuous block backward differentiation formula for solving stiff ordinary differential equations, Comput. Math. Appl. 65 (2013) pp.996–1005.

[22] S.G. Pinto, D.H. Abreu, B. Simeon, Strongly A-stable first stage explicit collocation methods with stepsize control for stiff and differential–algebraic equations, J. Comput. Appl. Math. 259 (2014) pp.138–152.

[23] X. Piao, S. Bu, D. Kim, P. Kim, An embedded formula of the Chebyshev collocation method for stiff problems, J. Comput. Phys. 351 (2017) pp.376-391.

[24] B. Wang, X. Wu, Exponential collocation methods for conservative or dissipative systems, J. Comput. Appl. Math. 360 (2019) pp.99–116.

[25] MM. Rashidi, G. Domairry, New analytical solution of the three-dimensional Navier-Stokes equations, Mod. Phys. Lett. B 23 (2009) pp. 3147-3155.

[26] H. Tunc, M. Sari, An efficient local transform method for initial value problems, Sigma Journal of Engineering and Natural Sciences 39(1) (2019) pp.163-174.

[27] H. Tunc, M. Sari, A local differential transform approach for the cubic nonlinear Duffing oscillator with damping term, Sci. Iran. 26(2) (2019) pp.879-886.

[28] I.H. Abdel-Halim Hassan, Differential transformation technique for solving higher-order initial value problems, Appl. Math. Comput. 154(2) (2004) pp.299-311.

[29] B. Liu, X. Zhou, Q. Du, Differential transform method for some delay differential equations, Applied Mathematics 6(3) (2015) pp.585-593.

[30] M.J. Jang, C.L. Chen, Y.C. Liy, On solving the initial-value problems using the differential transformation method, Appl. Math. Comput. 115(2-3) (2000) pp.145-160.

[31] L.F. Shampine, M.W. Reichelt, The MATLAB ODE Suite, Siam J Sci Comput 18 (1997) pp.1–22.

[32] J. Nocedal, S.J. Wright, Numerical Optimization, Springer, 2006.

[33] R. Li, S. Pei, B. Chen, Y. Song, T. Zhang, W. Yang, J. Shaman, Substantial undocumented infection facilitates the rapid dissemination of novel coronavirus (SARS-CoV-2), Science 368 (2020) pp.489–493.

[34] Z.M. Odibat, C. Bertelle, M.A. Aziz-Alaoui, G.H.E. Duchamp, A multi-step differential transform method and application to non-chaotic or chaotic systems, Comput. Math. Appl. 59(4) (2010) pp.1462-1472.

[35] E.R. El-Zahar, An adaptive step-size Taylor series based method and application to nonlinear biochemical reaction model, Trends Appl. Sci. Res. 7(11) (2012) pp.901-912.